\documentclass[centertags,12pt]{amsart}
\usepackage{latexsym}
\usepackage{amssymb}
\usepackage{amscd}
\numberwithin{equation}{section}
\makeatletter
\renewcommand{\subsection}{\@startsection
{subsection}{2}{0mm}{\baselineskip}{-0.25cm}
{\normalfont\normalsize\bf}}
\makeatother
\newtheorem{theorem}{Theorem}[section]
\newtheorem{proposition}[theorem]{Proposition}
\newtheorem{lemma}[theorem]{Lemma}
\newtheorem{corollary}[theorem]{Corollary}
{\theoremstyle{remark}
\newtheorem*{claim*}{Claim}
\newtheorem{remark}[theorem]{Remark}}
\theoremstyle{definition}
\newtheorem{example}[theorem]{Example}
\newtheorem*{question*}{Question}
\def\F{\mathbb F}
\def\P{\mathbb P}
\def\N{\mathbb N}

\def\cC{\mathcal C}
\def\cD{\mathcal D}
\def\cF{\mathcal F}
\def\cH{\mathcal H}

\def\cK{\mathcal K}

\def\cW{\mathcal W}
\def\cX{\mathcal X}
\def\cY{\mathcal Y}
\def\l{\ell}
\def\lq{\mathbb L_q}
\def\fq{\mathbb F_q}
\def\sq{\sqrt{q}}
\def\supp{{\rm Supp}}
\def\div{{\rm div}}
\def\dim{{\rm dim}}
\def\deg{{\rm deg}}
\def\frx{{\rm Fr}_{\mathcal X}}
\begin{document}
\author[Cossidente]{A. Cossidente}
\author[Korchm\'aros]{G. Korchm\'aros}
\author[Torres]{F. Torres}\thanks{1991 Math. Subj. Class.: Primary 11G,
Secondary 14G}\thanks{This research was
carried out within the activity of
GNSAGA of the Italian CNR with the support of the Italian
Ministry for Research and Technology. Torres acknowledges support of
Cnpq-Brazil}
\title[Curves covered by the Hermitian curve]{On curves covered by
the Hermitian curve}
\address{Dipartimento di Matematica 
Universit\`a della Basilicata 
via N. Sauro 85, 85100 Potenza, Italy}
\email{cossidente@unibas.it}
\email{korchmaros@unibas.it}
\address{IMECC-UNICAMP, Cx. P. 6065, Campinas-13083-970-SP, Brazil}
\email{ftorres@ime.unicamp.br}
     \begin{abstract}
For each proper divisor $d$ of $(q-\sq+1)$, $q$ being a square power of
a prime, maximal curves $\fq$-covered by the Hermitian curve of genus 
$\frac{1}{2}(\frac{q-\sq+1}{d}-1)$ are constructed. 
     \end{abstract}
\maketitle
\section{Introduction}\label{s1}
A maximal curve over a finite field $\fq$, is a 
projective geometrically irreducible non-singular algebraic curve   
defined over $\fq$ whose number of $\fq$-rational points attains the 
Hasse-Weil upper bound
$$
q+1+2g\sq\, ,
$$
where $g$ is the genus of the curve. These curves were studied in 
\cite{r-sti}, \cite{sti-x}, \cite{ft}, \cite{geer-vl1} (see also the 
references therein), \cite{fgt}, \cite{ft1}, \cite{geer-vl2}, \cite{gt}, 
and \cite{chkt}. Maximal curves have been intensively studied also in 
connection with Coding Theory in which such curves play an important role 
\cite{go}.

The Hermitian curve $\cH$, namely the plane curve defined by 
    \begin{equation}\label{eq1.1}
Y^{\sq}Z+YZ^{\sq}=X^{\sq+1}\, ,
    \end{equation}
is a well-known example of a maximal curve over $\fq$. By a result of
Lachaud \cite[Proposition 6]{lachaud}, any non-singular curve 
$\fq$-covered by $\cH$ is also maximal:  indeed almost all of the known 
maximal curves arise in this way. Lachaud's result has pointed out two
fundamental (and still open) problems, namely the classification problem
of maximal curves covered by $\cH$ and the existence problem of
maximal curves not covered by $\cH$.

In this paper, we present a new infinite class of maximal curves
$\fq$-covered by $\cH$, which arises from a cyclic automorphism group of
$\cH$ of order $(q-\sq+1)$, and contains a curve of genus 
$\frac{1}{2}(\frac{q-\sq+1}{d}-1)$ for each proper divisor $d$ of 
$(q-\sq+1)$. The construction together with some remarkable properties are
given in Section \ref{s5} and Section \ref{s6}. Since the $\fq$-invariant
linear series $\cD:=|(\sq+1)P_0|$ plays an important role in some current
research on maximal curves over $\fq$, we have computed its dimension and
$(\cD,P)$-orders for some points $P$. These curves are non-classical (for
the canonical morphism) whenever $d\le (q-2\sq+1)/(2\sq-1)$. Moreover, 
for small values of $d$, especially for $d=3$ (and 
$\sq\equiv\! 2\pmod{3}$) our results are also related to recent 
investigations on maximal curves having many rational points, or 
equivalently, large genus with respect to $q$. 
Such maximal curves are somewhat rare. However their classification is a
difficult task and still in progress. The first general result in this
direction states that the genus of a maximal curve over $\fq$ is at most 
$\sq(\sq-1)/2$ with equality holding if and only if the curve is
$\fq$-isomorphic to the Hermitian curve, see \cite{ihara}, \cite{r-sti}. 
In \cite{ft} it is proved that no maximal curve over $\fq$ has genus $g$
with $g \in ](\sq-1)^2/4,\sq(\sq-1)/2[$, a result conjectured in
\cite{sti-x}. From \cite[Theorem 3.1]{fgt}, the second largest genus of 
maximal curves over $\fq$ , $q$ odd, is equal to $(\sq-1)^2/4$, and such  
curves are $\fq$-isomorphic to the Artin-Schreier curve 
$y^{\sq}+y=x^{(\sq+1)/2}$. For $q$ even the second largest genus 
is $\sq(\sq-2)/4$, and it seems that 
curves having such a genus are also unique up to $\fq$-isomorphism (cf. 
\cite{at}). 
Instead, the problem of determining the third large genus seems to be much
more involved and heavily dependent on certain arithmetical behaviour of 
$q$  as the list after Remark \ref{remark6.1} suggests. From 
\cite[Proposition 2.5]{ft1}, the third large genus , for $q$ odd, is at
most $(\sq-1)(\sq-2)/4$. Good candidates come from the known infinite classes
of maximal curves $\fq$-covered by $\cH$, namely $\cF_2$ and $\cH_3$, 
where $\cF_t$ is the Fermat curve $x^{(\sq+1)/t}+y^{(\sq+1)/t}+1=0$ and 
$\cH_t$ the Artin-Schreier curve $y^{\sq}+y=x^{(\sq+1)/t}$, $t$ being a 
proper divisor of $(\sq+1)$. It is worth mentioning that both $\cF_2$ and
$\cH_4$ have genus $(\sq-1)(\sq-3)/8$, but they are not $\fq$-isomorphic,  
\cite{chkt}. Also the curve $\cH_t$ has been 
characterized via the type of Weierstrass semigroups at $\fq$-rational
points \cite[Theorem 2.3]{fgt}. We conjecture that the candidate for the
third largest genus in the case of $\sq\equiv 2\!\pmod{3}$ is the maximal
curve of genus $(q-\sq-2)/6$ described in Section \ref{s6}. This emerges
from Section \ref{s3} where maximal curves with $\dim(\cD)=3$ are 
investigated.

As in previous papers concerning maximal curves \cite{ft}, \cite{fgt},
\cite{ft1}, \cite{gt}, and \cite{chkt} we have used St\"ohr-Voloch's
approach to the Hasse-Weil bound to carry on our research. 
\smallskip

\noindent {\bf Conventions.} The word curve means a projective 
geometrically irreducible algebraic curve defined over a finite field 
$\fq$. For a curve $\cX$, $\bar\cX$ denotes its nonsingular model over 
$\fq$. 
\section{Background}\label{s2}
In this section we collect some results concerning Weierstrass Point
Theory, Frobenius orders and maximal curves. 
\subsection{Weierstrass Point Theory and Frobenius orders}\label{s2.1}
Here we recall relevant material from St\"ohr-Voloch's \cite[\S1,\S2]{sv}. 

Let $\cX$ be a nonsingular curve of genus $g$ defined over $\bar\fq$
equipped with the action of the Frobenius morphism $\frx$ over $\fq$, and
let $\cD$ be a $g^r_d$ on $\cX$. Suppose that $\cD$ is defined over $\fq$.

Associated to $\cD$ we have two divisors on $\cX$, namely the {\em
ramification divisor} $R=R^{\cD}$, and the {\em $\fq$-Frobenius divisor}
$S=S^{(\cD,q)}$. These divisors provide a lot of geometrical and
arithmetical information on $\cX$. We recall that the set of Weierstrass
points $\cW_{\cX}$ of $\cX$ is the support of $R^{\cK}$, where
$\cK=\cK_{\cX}$ is the canonical linear series on $\cX$. 

For $P\in \cX$ let us denote by $j_i(P)=j_i^{\cD}(P)$ the $i$-th
$(\cD,P)$-order, by $\epsilon_i=\epsilon_i^{\cD}$ the $i$-th $\cD$-order
($i=0,\ldots,r$), and by $\nu_i=nu_i^{(\cD,q)}$ the $i$-th $\fq$-Frobenius
order of $\cD$ ($i=0,\ldots,r-1$). The following are the main properties
of $R$ and $S$. Set $p:={\rm char}(\fq)$.
\begin{enumerate}
\item\label{1} $\deg(R)=(2g-2)\sum_{i=0}^{r}\epsilon_i + (r+1)d$;
\item\label{2} $j_i(P)\ge \epsilon_i$ for each $i$ and each $P$;
\item\label{3} $v_P(R)\ge \sum_{i=0}^{r} (j_i(P)-\epsilon_i)$ and equality
holds iff ${\rm det}(\binom{j_i(P)}{\epsilon_j}) \not\equiv 0\pmod{p}$;
\item\label{4} $(\nu_i)$ is a subsequence of $(\epsilon_i)$;
\item\label{5} $\deg(S)= (2g-2)\sum_{i=0}^{r-1}\nu_i + (q+r)d$;
\item\label{6} For each $i$ and for each $P\in \cX(\fq)$, $\nu_i\le
j_{i+1}(P)-j_1(P)$;
\item\label{7} For each $P\in \cX(\fq)$, $v_P(S)\ge
\sum_{i=0}^{r-1}(j_{i+1}(P)-\nu_i)$ and equality holds iff ${\rm
det}(\binom{j_{i+1}(P)}{\nu_j})\not\equiv 0 \pmod{p}$.
\end{enumerate}
Therefore if $P\in \cX(\fq)$, \ref{6}) and \ref{7}) imply
\begin{enumerate}
\item[8.] $v_P(S)\ge rj_1(P)$.
\end{enumerate}
Consequently from \ref{5}) and 8) we obtain the main result of 
\cite{sv}, namely 
\begin{enumerate}
\item[9.] $\# \cX(\fq)\le \deg(S)/r$.
\end{enumerate}
\subsection{Maximal curves}\label{s2.2} We summarize some results 
from \cite{fgt} and \cite{ft1}. Let $\cX$ be a maximal curve over $\fq$. 
The key property on $\cX$ is the following linear equivalence 
\cite[Cor. 1.2]{fgt} 
\begin{equation}\label{eq2.1}
\sq P+\frx(P)\sim (\sq+1)P_0\, ,\quad\qquad \text{$P\in \cX$,\ $P_0\in 
\cX(\fq)$}\, .
\end{equation}
Hence, for $P_0\in \cX(\fq)$, $\cX$ is equipped with the $\fq$-invariant  
linear series 
$$
\cD_{\cX}:=|(\sq+1)P_0|\, ,
$$
so that $\dim(\cD_{\cX})$ is independent of  
$P_0\in X(\fq)$. We have that $\dim(\cD_{\cX})\ge 2$ \cite[Prop. 
1]{sti-x} (see  
also \cite[Prop. 1.5(iv)]{fgt}). Furthermore we have the
        \begin{lemma}\label{lemma2.1} Let $\cX$ be a maximal
curve of genus $g$ over $\fq$. The following statements are 
equivalent 
\begin{enumerate}
\item $\cX$ is $\fq$-isomorphic to the Hermitian curve (so that
$g=\sq(\sq-1)/2$);
\item $g>(\sq-1)^2/4$;
\item $\dim(\cD_{\cX})=2$.
\end{enumerate}
        \end{lemma}
\begin{proof} See \cite{r-sti} and \cite[Thm. 2.4]{ft1}.
\end{proof}
Set $n+1:=\dim(\cD_{\cX})$, $j_i(P):=j^{\cD_\cX}_i(P)$,
$\epsilon_i:=\epsilon_i^{\cD_\cX}$, $\nu_i:=\nu_i^{(\cD_\cX,q)}$ and
denote by $m_i(P)$ the $i$-th non-gap at $P\in \cX$. 
So from (\ref{eq2.1}) and \cite[\S2.3]{ft1}, for each $P\in \cX$ the
following holds:
\begin{equation}\label{eq2.2}
m_0(P)=0<m_1(P)<\ldots<m_n(P)=\sq\, .
\end{equation}
Moreover, $m_{n+1}(P)=\sq+1$ if $P\in \cX(\fq)$, and $m_{n-1}=\sq-1$ if
$P\not\in \cW_{\cX}$ and $n+1\ge 3$ \cite[Prop. 1.5(iv)(v)]{fgt}. The
main properties of the $(\cD_\cX,P)$-orders, $\cD_\cX$-orders and
$\fq$-Frobenius orders of $\cD_\cX$ are the following (see
\cite[\S1]{fgt}, \cite[\S2.2]{ft1}).
\begin{enumerate} 
\item For each $P\in \cX$, $j_1(P)=1$ so that $\epsilon_1=1$;
\item $\epsilon_{n+1}=\nu_n=\sq$;
\item If $P\in \cX(\fq)$, then the $(\cD_\cX,P)$-orders are 
$\sq+1-m_i(P)$ ($i=0,1,\ldots,n+1$) so that $v_P(R^{\cD_\cX}\ge1$;
\item If $P\not\in \cX(\fq)$, then $j_{n+1}(P)=\sq$ and there exists 
$I=I(P)\in [1,n]$ such that
$$
\sq-m_n(P)<\ldots<\sq-m_{n-I+1}<j_I(P)<\sq-m_{n-I}<\ldots<\sq-m_0(P)
$$
are the $(\cD_\cX,P)$-orders.
\item In particular if $P\not\in \cW_{\cX}$, then $\tilde m_i=m_i(P)$ is
independent of $P$; hence $\sq-\tilde m_i$ is a $\cD$-order for
$i=0,1,\ldots,n$. 
\item $\nu_i=\sq-\tilde m_{n-i}$ for $i=0,\ldots,n$, so that $\nu_1=1$
whenever $n+1\ge 3$.
\end{enumerate}
\section{On maximal curves with $\dim(\cD_{\cX})=3$}\label{s3}
Let $\cX$ be a maximal curve over $\fq$ of genus $g$.  We keep the
notations of the previous section. To study $\cX$, by Lemma
\ref{lemma2.1}, we can assume that $g\le (\sq-1)^2/4$ or equivalently that
$\dim(\cD_\cX)\ge 3$. This section deals with the case $\dim(\cD_\cX)=3$. 
So the $\cD_\cX$-orders are $\epsilon_0=0, \epsilon_1=1, \epsilon_2$ and
$\epsilon_3=\sq$, and the $\fq$-Frobenius orders of $\cD$ are $\nu_0=0$,
$\nu_1=1$ and $\nu_2=\sq$ (cf. \S\ref{s2.2}(1)(5)(6)). 

We first state a sufficient condition to have $\dim(\cD_\cX)=3$ and to 
compute $\epsilon_2$.
\begin{lemma}\label{lemma3.1} Let $X$ be a maximal curve over $\fq$ of
genus $g$ such that 
$$
(\sq-1)(\sq-2)/6<g\le (\sq-1)^2/4\, .
$$ 
Then 
\begin{enumerate} 
\item $\dim(\cD_\cX)=3$.  
\item $\epsilon_2\le 3$ and $\epsilon_2=2$ provided that $p={\rm
char}(\fq)\neq 3$.  
\end{enumerate} 
\end{lemma}
\begin{proof} (1) Since $\sq, \sq+1\in H(P)$ for each $P\in \cX(\fq)$,
(see \S\ref{s2.2}), $\cD_\cX$ is 
simple and hence we can apply Castelnuovo's genus bound for curves in
projective spaces as given in \cite[p.34]{fgt}: $g$ satisfies 
$$
2g\le \begin{cases} 
\frac{(q-n/2)^2}{n} & \text{if $n$ is even}\,,\\ 
\frac{(q-n/2)^2-1/4}{n} & \text{otherwise}\, , 
\end{cases}
$$
where $n+1=\dim(\cD_\cX)$. Therefore $n+1\ge 4$ would imply $g\le
(q-1)(q-2)/6$, a contradiction. Then we have $n+1\le 3$ and Lemma
\ref{lemma2.1} implies $n+1=3$.

(2) Let $P\in \cX(\fq)$. By \S\ref{s2.1}(7)(2), $v_P(S)\ge
j_2(P)+1\ge \epsilon_2+1$. This inequality, the maximality of $\cX$ 
and \S\ref{s2.1}(5) imply
\begin{equation*}
(\sq+1)(2g-2)+(q+3)(\sq+1)\ge 
(\epsilon_2+1)[(\sq+1)^2+\sq(2g-2)]\, .\tag{$*$}
\end{equation*}
If $\epsilon_2\ge 4$, then we would have 
$$
(\sq+1)(q-5\sq-2)\ge (2g-2)(4\sq-1)\, ,
$$
and from the upper bound on $g$ it follows that $0>q-2\sq+10$, a
contradiction.

Now $\epsilon_2=3$ implies $\binom{3}{2}\equiv 0\pmod{p}$ by the 
$p$-adic criterion \cite[Cor. 1.9]{sv} and hence $p=3$. This finishes the
proof of (2).  
\end{proof}
\begin{remark}\label{remark3.0} Here we show that the hypothesis on the
genus in Lemma \ref{lemma3.1}(1) is sharp whenever $\sq\equiv 0\pmod{3}$.

Let $q=p^m$ with $m$ even and $1\le r\le
m/2$. van der Geer and van der Vlugt \cite[Thm. 3.1]{geer-vl1}, 
\cite[Remark 5.2]{geer-vl2} constructed a maximal curve $\cX$ over $\fq$ 
of genus $(p^r-1)\sq/2$, by considering fibre products of curves of
type $y^p-y=ax^{\sq+1}$ with $a\in \fq^*$ satisfying $a^{\sq}+a=0$. It is
not difficult to see that a plane model for $\cX$ is given by an equation
of type 
$$
\sum_{i=0}^{r}y_1^{p^i} = b x^{\sq+1}\, ,
$$
where $b\in \fq^*$ and $b^{\sq}+b=0$ so that $\cX$ is $\fq$-covered by the
Hermitian curve. (The case $r=m/2$ has worked out in \cite[\S V, Example
E]{g-sti}). 

Set $\cD:=\cD_{\cX}$ and let $P_0\in \cX(\fq)$ be the unique point
over $x=\infty$. Now we can apply the case $``mn=\ell"$ in \cite{ft1}
to show the following facts:
    \begin{enumerate}
\item \quad $\dim(\cD)=p^{m/2-r}+1$;
\item The $(\cD,P)$-orders are as follows:

(i) 0, $\sq+1-ip^r$, $i=0,1,\ldots, p^{m/2-r}$ if $P=P_0$;

(ii) $0,1,\ldots, p^{m/2-r}$ and $\sq+1$ if $P\in
\cX(\fq)\setminus\{P_0\}$; 

(iii) $0,1,\ldots,p^{m/2-r}$ and $\sq$ if $P\in \cX\setminus\cX(\fq)$.
    \end{enumerate}
In particular, (putting $p=3$ and $r=m/2-1$) there exists a maximal curve
$\cX$ of genus $\sq(\sq-3)/6$ such that $\dim(\cD_{\cX})=4$. This shows
the sharpness of the upper bound for the genus in Lemma \ref{lemma3.1}(1)
for $\sq\equiv 0\pmod{3}$. 
\end{remark}
\begin{remark}\label{remark3.1} Let $\cX$ be a maximal curve over $\fq$ of
genus $g$ such that $\dim(\cD_{\cX})=3$.

(1) If $\epsilon_2$=2, then 
\S\ref{s2.1}(1), \S\ref{s2.2}(3) and the maximality of $\cX$ imply $g\ge 
(q-2\sq+3)/6$.

(2) If $\epsilon_2=3$, then relation $(*)$ in the proof of Lemma
\ref{lemma3.1} implies\\ 
$2g-2\le (\sq+1)(q-4\sq-1)/(3\sq-1)$. Notice that the upper bound on $g$
in Lemma \ref{lemma3.1} is equivalent to $2g-2\le (\sq+1)(\sq-3)/2$. 
\end{remark}
Next we investigate the values of $m_1(P)$ for $P\in \cX(\fq)$ when
$\dim(\cD_\cX)=3$. Notice that $m_1(P)\ge \sq/2$ since $2m_1(P)\in H(P)$
and $m_2(P)=\sq$, and that $g$ is bounded from above by the genus of the
semigroup $\langle m_1(P),\sq,\sq+1\rangle$. We have the following result
due to Fuhrmann (see also \cite[\S3.II]{selmer}). 
\begin{lemma}[F,\S A2]\label{lemma3.2} Let $\l,m\in \N$ such that $m/2\le
\l<m$. 
Let $\tilde g$ be the genus of the semigroup $\langle \l, m,m+1\rangle$. 
\begin{enumerate}
\item If $\l\not\in\{\lfloor\mbox{$\frac{m+1}{2}$}\rfloor,
m-1,\lfloor\mbox{$\frac{2m+2}{3}$}\rfloor,m-2\}$, then
$$
\tilde g\le \begin{cases}
\max{((m^2+4)/8,(m^2+3m)/10)} & \text{if $\l\le 3m/5$}\, ,\\
(m^2-5m+24)/6 & \text{if $3m/5\le 
\l<\lfloor\mbox{$\frac{2m+2}{3}$}\rfloor$}\, ,\\
(m^2-7m+70)/6 & \text{if $m\equiv 2 \pmod{3}$ and}\\ 
              & (2m+5)/3\le \l\le m+1-\sqrt{m+1}\, ,\\
(m^2-5m+40)/6 & \text{if $m\equiv 1 \pmod{3}$ and}\\
              & (2m+4)/3\le \l\le m+1-\sqrt{m+1}\, ,\\
(m^2-3m+18)/6 & \text{if $\sq\equiv 0\pmod{3}$ and}\\ 
              & (2m+3)/3\le \l\le m+1-\sqrt{m+1}\, ,\\
(m^2+2m+9)/8  & \text{if $m+1-\sqrt{m+1}<\l<m-2$}\, .
\end{cases}
$$
\item If $\l\in\{\lfloor\mbox{$\frac{m+1}{2}$}\rfloor,m-1\}$, then
$\tilde g=(m-1)^2/4$.
\item If $\l\in\{\lfloor\mbox{$\frac{2m+2}{3}$}\rfloor,m-2$, then
$$
\tilde g=
\begin{cases}
(m^2-m+4)/6 & \text{for $m\equiv 2\pmod{3}$}\, ,\\
(m^2-m)/6 & \text{otherwise}\, .
\end{cases}
$$
\end{enumerate}
\end{lemma}
        \begin{corollary}\label{cor3.1} Let $\cX$ be maximal curve over 
$\fq$. Suppose that $\dim(\cD_\cX)=3$ and that $\epsilon_2=2$. Then for 
$P\in \cX(\fq)$ we have that
$$
m_1(P)\in \{\lfloor\mbox{$\frac{\sq+1}{2}$}\rfloor,
\sq-1,\lfloor\mbox{$\frac{2\sq+2}{3}$}\rfloor,\sq-2\}\, .
$$
         \end{corollary}
\begin{proof} The genus $g$ of $\cX$ is bounded from above by the genus of
the semigroup $\langle m_1(P),\sq,\sq+1\rangle$. Then the result follows
by applying the lemma with $m=\sq$ taking into consideration 
that  $g\ge (q-2\sq+3)/6$ (Remark \ref{remark3.1}(1)).
\end{proof}
\begin{remark}\label{remark3.2} Let $\cX$ be a maximal curve over $\fq$
and let $P\in \cX(\fq)$.

(1) If $q$ is odd and $m_1(P)=(\sq+1)/2$, then $\dim(\cD_\cX)=3$,
$\epsilon_2=2$ and $\cX$ is $\fq$-isomorphic to the nonsingular model of
$y^{\sq}+y=x^{\sq+1}$ \cite[Thm. 2.3]{fgt}.

(2) If $q$ is even and $m_1(P)=\sq/2$, then $\dim(\cD_\cX)=3$,
$\epsilon_2=2$ and $\cX$ is $\fq$-isomorphic to the nonsingular model of
$\sum_{i=0}^{t}y^{\sq/2^i}=x^{\sq+1}$, $\sq=2^t$ \cite{at}.
\end{remark}
         \begin{lemma}\label{lemma3.3} 
Let $\cX$ be a maximal curve over $\fq$ such that $\dim(\cD_{\cX})=3$.
Then there exists $P\in \cX(\fq)$ such that $m_1(P)+\epsilon_2=\sq+1$.
         \end{lemma} 
         \begin{proof}
This follows from the proof of \cite[Prop.  1.5(v)]{fgt}; for the sake of
completeness we write a proof. By \S\ref{s2.2}(3), for each $P\in
\cX(\fq)$, $j_2(P)=\sq+1-m_1(P)$ so it will be enough to show that
$j_2(P)=\epsilon_2$ for some $\fq$-rational point. Suppose that
$j_2(P)>\epsilon_2$ for each $P\in \cX(\fq)$ so that $v_P(R^{\cD_\cX})\ge
2$ (see \S\ref{s2.1}(2)(3)). Using $\epsilon_2\le \sq-1$ and the
maximality of $\cX$ from \S\ref{s2.1}(1) we have 
$$ 
2\sq(2g-2)+4(\sq+1)\ge 2\sq(2g-2)+2(\sq+1)^2\, , 
$$ 
i.e. $0\ge (\sq+1)(2\sq-2)$, a contradiction.
\end{proof}
It follows immediately the
          \begin{corollary}\label{cor3.2} 
Let $\cX$ be a maximal curve over $\fq$ with $\dim(\cD_\cX)=3$. Then the
number $(\sq-2)$ (resp.  $(\sq-1)$) is realized as a non-gap at a
$\fq$-rational point of $\cX$ iff $\epsilon_2=3$ (resp.  $\epsilon_2=2$). 
          \end{corollary}
          \begin{remark}\label{remark3.3}
From Remark \ref{remark3.2}, the number $\sq-1$ is always realized as a
non-gap at $\fq$-rational points of maximal curves. So far we do not know
any example of a maximal curve over $\fq$ with $\dim(\cD_\cX)=3$ and
having a $\fq$-rational point $P$ with $m_1(P)\in
\{\lfloor\mbox{$\frac{2\sq+2}{3}$}\rfloor, \sq-2\}$. 
          \end{remark}
          \begin{remark}\label{remark3.4}
If $q$ is a square and $\sq\equiv 0 \pmod{3}$, there exists a maximal
curve $\cX$ over $\fq$ whose genus is $\sq(\sq-1)/6$. Indeed this is a
particular case of the maximal curves constructed in 
\cite[Prop. 5.1]{geer-vl2}. From Lemma \ref{lemma3.1} it follows that 
$\dim_{\cD_\cX}=3$. 
           \end{remark}
We finish this section giving a characterization of the set
$\supp(S^{\cX})$. Let $n+1=\dim(\cD_\cX)$, $\cX$ being a maximal curve
over $\fq$. Let $P\in\cX\setminus\cX(\fq)$. From the proof of \cite[Thm.
2.1]{gt} we see that $P\in \supp(S^{\cD_\cX})\Rightarrow 
m_1(P)<\sq-n+1$. As an scholium of (loc. cit.), for $n+1=3$ 
the converse also holds:
          \begin{lemma}
Let $\cX$ be a maximal curve over $\fq$ with $\dim(\cX)=3$. Then
$$
\supp(S^{\cD_\cX})=\cX(\fq)\cup\{P\in\cX\setminus \cX(\fq):
m_1(P)<\sq-1\}\, .
$$
          \end{lemma}
          \begin{proof} 
We already know that $\cX(\fq)\cup\{P\in \cX\setminus
\cX(\fq):m_1(P)<\sq-1\}\subseteq \supp(S^{\cD_\cX})$ by \S\ref{s2.1}(8)
and the remark stated before the lemma. Conversely let $P\in\cX\setminus
\cX(\fq)$ with $m_1(P)<\sq-1$. Then the $(\cD,P)$-orders are $j_0=0, j_1= 
1,j_2=q-m_1(P)$ and $\sq$ (see \S\ref{s2.2}(4)). Let $u, v, w\in 
\bar\fq(\cX)$ such that $j_1P+D_1=\div(u)+(\sq+1)P_0$ with
$P\not\in\supp(D_1)$; 
$\div(v)=D_v-m_1(P)P$ with $P\not\in\supp(D_v)$  
and $\div(w)=\sq P+\frx(P)-(\sq+1)P_0$ (cf. (\ref{eq2.1}). 
Then
$$
\div(uv)+(\sq+1)P_0=D_u+j_2(P)P+\frx(P)\ \text{and}\  
\div(v)+(\sq+1)P_0=\sq P +\frx(P)
$$
so, according to the proof of \cite[Thm. 1.1]{sv}, $\frx(P)$ belongs to
the tangent line at $P$. This line is generated by $\pi(P)$ and
$(D^1_t\pi)(P)$, where $\pi=(1:u:uv:w)$,
$D^1_t\pi=(0:D^1_tu:D^1_t(uv):D^1_tw)$ and $t$ is a local parameter at
$P$. Since the $\fq$-orders are $0,1$ and $\sq$ (see \S\ref{s2.2}(6)) it 
follows that $P\in \supp(S^{\cD_\cX})$.
          \end{proof}
\section{Plane models for the Hermitian curve}\label{s4} 
The aim of this section is to introduce two further plane models for the
Hermitian curve $\cH$ defined by (\ref{eq1.1}). 
\subsection{A plane singular model for $p:={\rm char}(\fq)\ge 
3$}\label{s4.1} 
The Hermitian curve has a cyclic automorphism $\psi$ of order $(q-\sq+1)$.
Using 
(\ref{eq1.1}), this automorphism becomes a linear collineation. 
Unfortunately, the associated matrix is not in diagonal form, and the
model (\ref{eq1.1}) is not appropriate to investigate those properties of
the Hermitian curve which depend on $\psi$. This was recognized at first
in \cite{c}, \cite{ck1}, \cite{ck} in the study of $\psi$-invariant arcs,
or, equivalently, of arcs that are the complete intersection of two
Hermitian curves in a suitable mutual position. In \cite{ck1}, a useful
singular plane model was considered for the Hermitian curve. Actually,
this model arises from the algebraic envelope of a $\psi$-invariant arc,
via Segre's fundamental theorem of $k$-arcs, see \cite{segre}, \cite[Thm.
10.4.3]{h}. The same model allows us to determine the quotient curve with
respect to the automorphism $\psi^{(q-\sq+1)/d}$, for each proper divisor
$d$ of $(q-\sq+1)$. These curves are maximal and are investigated in 
\S\ref{s5}. Here we limit ourselves to describe how an equation for this
model can be obtained. The starting point is
the following lemma which follows from the main result in \cite{ck1} and
\cite{r-sti}.  
         \begin{lemma}\label{lemma4.1.1} 
The nonsingular model $\bar\cC$ of the algebraic envelope $\cC$ of a 
$\psi$-invariant arc in $\P^2(\fq)$, regarded as a curve of degree
$2(\sq+1)$ in the dual plane of $\fq$, is $\fq$-isomorphic to the
Hermitian curve over $\fq$.
         \end{lemma}
Let
$$
\alpha: \P^2(\bar\fq)\to \P^2(\bar\fq)\, ,\quad (x:y:z)\mapsto
(ax:a^{q+1}y:z)\, ,
$$
where $a\in \F_{q^3}$ is a primitive $(q^2+q+1)$-th root of unity, and let
$$
E:=(1:1:1)\, .
$$ 
Then the orbit of $E$ under $\alpha$ is given by
$$
\Pi=\{(a^i:a^{(q+1)i}:1): i=0,1,\ldots,q^2+q\}\, .
$$ 
        \begin{lemma}[CK, Prop. 1]\label{lemma4.1.2} 
The set $\Pi$ is a projective subplane of $\P^2(\F_{q^3})$ lying in a
non-classical position, i.e.  $\Pi\neq \P^2(\fq)$. 
        \end{lemma}
Then $\Pi=\P(\lq^2)$ with $\lq$ a field isomorphic to $\fq$. Now if $\cC:
F(X,Y,Z)=0$, with $F\in\fq[X,Y,Z]$, one writes the corresponding curve
$\cC'$ over
$\bar\lq=\bar\F_{q^3}$ as follows. Let 
$$
A_1:=(1:0:0)\, ,\quad A_2:=(0:1:0)\, ,\quad A_3:= (0:0:1)\, ,
$$ 
and choose three points $A_1', A_2', A_3'\in \Pi$ such that
$A_1'A_2'A_3'E$ is a non-degenerate quadrangle. Let
$$
\kappa:\P^2(\bar\fq)\to \P^2(\bar\F_{q^3})
$$ 
be the projective map such that $\kappa(A_i)=A_i'$ ($i=1,2,3$) and
$\kappa(E)=E$. Then $\cC':= \kappa(\cC)$. Let
\begin{align*} 
\beta &: \P^2(\bar\fq)\to \P^2(\bar\fq)\, ,
\quad (x:y:z)\mapsto (z:x:y)\, , \\
\intertext{and}
\gamma &: \P^2(\bar\fq)\to \P^2(\bar\fq)\, ,
\quad (x:y:z)\mapsto (x^q:y^q:z^q)\, .
\end{align*}
        \begin{lemma}\label{4.1.3}
\begin{enumerate}
\item $\Pi$ is the set of fixed points of $\beta\circ\gamma$, i.e. 
$\beta\circ\gamma$ is the Frobenius morphism in the 
new frame $A_1'A_2'A_3'E$. 
\item The curve $\cC'$ is defined over $\lq$ iff $\beta(\cC)=\cC$.
\end{enumerate}
        \end{lemma}
        \begin{proof} (1) It is clear that $\Pi$ is contained in the 
set of fixed points of $\beta\circ\gamma$. Conversely let $P=(b:c:1)$ be 
a fixed point of $\beta\circ\gamma$; then $(1:b^q:c^q)=(b:c:1)$ so that 
$bc^q=1$ and $b^{q+1}=c$. Then $b^{q^2+b+1}=1$ and (1) follows.

(2) Let $G=G(X',Y',Z')=0$ be the equation of
$\cC'$, i.e. $F=G\circ\kappa\ (*)$. Then $\cC'$ is defined over $\lq$ iff 
$(G\circ\kappa)\circ(\beta\circ\gamma)={\rm Fr}\circ G\circ\kappa$ where 
${\rm Fr}$ is the Frobenius morphism on $\P^1(\bar\lq)$. So (2)  
follows from $(*)$ and the fact that $\cC$ is defined over $\fq$.
         \end{proof}
        \begin{lemma}\label{lemma4.1.3}
\begin{enumerate}
\item The curve $\cC'$ in coordinates $(X_0,X_1,X_2)=k(X,Y,Z)$ is defined
by 
\begin{align}\label{eq4.1}
\begin{split}
G(X_0,X_1,X_2)= & X_1^2X_2^{2\sq}+X_0^{2}X_1^{2\sq}+
X_0^{2\sq}X_2^{2}-\\
             & 2(X_0^{\sq+1}X_1^{\sq}X_2+
X_0^{\sq}X_1X_2^{\sq+1}+X_0X_1^{\sq+1}X_2^{\sq})=0\, .\\
\end{split}
\end{align}
\item The curve $\cC'$ is defined over $\lq$. 
\end{enumerate}
        \end{lemma}
        \begin{proof} (1) See \cite[Prop. 6]{ck1}.

(2) We know that $\cC$ is defined over $\fq$ and that $F=G\circ\kappa$. 
Then to apply Lemma 
\ref{lemma4.1.3} it is enough to show that
$G\circ\kappa\circ\beta=G\circ\kappa$. This easily can be checked using
(\ref{eq4.1}).
        \end{proof}
Starting from (\ref{eq4.1}) we can write an equation for $\cC$ over $\fq$
as follows. Let 
$$
A_1'=(a:a:^{q+1}:1)\qquad A_2'=(1:a:a^{q+1})\qquad A_3'=(a^{q+1}:1:a)\, ,
$$
where $a\in \F_{q^3}$ is as above. Since
$a^{-q-1}=a^{q^2}$ and $a^{-1}=a^{q^2+q}$, then $A_2', A_3'\in \Pi$. 
Let $\kappa$ be the projective map induced by the matrix
\begin{equation}\label{eq4.2}
M= \begin{bmatrix} 
a & 1 & a^{q+1} \\
a^{q+1} & a & 1 \\
1 & a^{q+1} & a
\end{bmatrix}\, .
\end{equation}
This map sends $A_i$ to $A_i'$ for $i=1,2,3$ and fixes $E$, and since 
$(a+1){\rm det}(M)=a^2+a+1$, $M$ is non-singular. 
Hence we have the main result of this subsection:
    \begin{proposition}\label{prop4.1.1}
A plane model over $\fq$ for the Hermitian curve $\cH$ is given by 
$$
H(X,Y,Z)=G(aX+Y+a^{q+1}Z, a^{q+1}X+aY+Z,X+a^{q+1}Y+aZ)=0\, ,
$$
with $G$ defined in (\ref{eq4.1}).
    \end{proposition}
\subsection{A plane non-singular model over $\F_{\protect\sq^3}$}\label{s4.2} 
It is well known that the Hermitian curve $\cH$ is projectively equivalent
over $\fq$ to the curve of equation
$$
H(X,Y,Z)=X^{\sq+1}+Y^{\sq+1}+Z^{\sq+1}\, .
$$
     \begin{proposition}\label{prop4.2.1} The Hermitian curve $\cH$ 
is $\F_{\sq^3}$-isomorphic to the projective plane curve $\cC'$
defined by 
$$ 
G(X_0,X_1,X_2)=X_0^{\sq} X_2 +X_2^{\sq} X_1+ X_1^{\sq} X_0\, .  
$$ 
     \end{proposition} 
We first prove the 
     \begin{lemma}\label{lemma4.2.1} 
There exists $a\in \F_{\sq^3}$ satisfying the following properties:  
\begin{enumerate} 
\item $a_1:=a^{q\sq+\sq}+a^{q+\sq+1}+a=0$;  
\item $a_2:=a^{q\sq+q+\sq+1}+a^{\sq+1}+1=0$;  
\item $a_3:=a^{q\sq+\sq+1}+a^{q+1}+a^{\sq}\neq 0$;  
\item The matrix $M$ in (\ref{eq4.2}) is nonsingular.
\end{enumerate}
     \end{lemma}
\begin{proof} We claim that a root $a$ of the polynomial 
$f(X)=X^{\sq+1}+X+1$ satisfies the lemma. In fact, from $a^{\sq+1}+a+1=0\ 
(*)$ we have $a^{q+\sq}+a^{\sq}+1=0\ (*1)$ so that
$a^{q+\sq+1}+a^{\sq+1}+a=0\ (*2)$. From $(*)$ and $(*2)$, $a^{q+\sq+1}=1\ 
(*3)$ so that $a^{\sq^3}=a\ (*4)$, i.e. $a\in \F_{\sq^3}$.

Now $a_1=a^{q\sq+\sq}+a^{q+\sq+1}+a=a^{1+\sq}+a^{q+\sq+1}+a=0$ by
$(*4)$ and $(*1)$;

$a_2=a^{q\sq+q+\sq+1}+a^{\sq+1}+1=a+a^{\sq+1}+1=0$, by $(*4)$, $(*3)$ and
$(*)$.

To prove (3), from  
$(*4), (*)$ and $(*3)$ we have $a(a+1)a_3=-(a^2+a+1)^2$ and so $a_3=0$ iff
$a^2+a+1=0$. Now $f(X)$ has $\sq+1$ different roots; then, as $\sq+1\ge
3$, we can pick a root $a$ of $f(X)$ such that $a^2+a+1\neq 0$.

Statement (4) follows from the identity $(a+1)^3{\rm
det}(M) = (a^2+a+1)^3$ and the proof of (3).
\end{proof}
\begin{proof} ({\it Proposition \ref{prop4.2.1}}) Let $a\in \F_{\sq^3}$ be
as in Lemma \ref{lemma4.2.1}. We use the notations of the preceding lemma.
Let 
$\kappa: \P^2(\bar\fq)\to \P^2(\bar \F_{\sq^3})$ be the projective  
map defined by $M$ in (\ref{eq4.2}). Then $\kappa^{-1}(\cC')$ is the
projective plane curve defined by 
$$
G(aX+Y+a^{q+1}Z, a^{q+1}X+aY+Z, X+a^{q+1}Y+aZ)=0\, .
$$
After some computations we obtain
\begin{align*}
G(aX+Y+a^{q+1}Z, a^{q+1}X+aY+Z, & X+a^{q+1}Y+aZ) =  a_3H(X,Y,Z)+\\
                                 & a_2H(X,Y,Z)+ a_1H(Y,X,Z)
\end{align*}
and we are done.
\end{proof}
\section{The case $d$ divides $(q$-$\protect\sq+1)$}\label{s5}
Throughout this section we assume that $q$ is a square power of an odd
prime and that $d$ is a positive divisor of $(q-\sq+1)$. The main result
here is the
     \begin{theorem}\label{thm5.1} Let $q$ be a square power of an odd
prime. Then for each positive divisor $d$ of $(q-\sq+1)$ there exists a
maximal curve over $\fq$ of genus $\frac{1}{2}(\frac{q-\sq+1}{d}-1)$.
     \end{theorem}
     \begin{remark}\label{remark5.1} If $\sq\equiv r 
\pmod{d}$ and $r^2-r+1\equiv 0 \pmod{d}$, then $d$ is odd and 
$\gcd(r,d)=1$. Then $d=3$ iff $r=2$ and if $d>3$, then $6$ divides
$\phi(d)$, $\phi$ being the Euler function. In particular, if $d$ is 
prime, then $d\equiv 1\pmod{6}$. 
     \end{remark}
To prove the theorem we use the model $\cC'$ over $\lq$ of the Hermitian
curve $\cH$ stated in Lemma \ref{lemma4.1.3}. We use the notations of 
subsection \ref{s4.1} and set  
$x:=X_0/X_2$, $y:=X_1/X_2$. Then, according to (\ref{eq4.1}), $\cC'$ is
defined by the affine equation
$$
y^2+x^2y^{2\sq}+x^{2\sq}-2(x^{\sq+1}y^{\sq}+x^{\sq}y +xy^{\sq+1})=0\, .
$$
We recall that $\cC'$ has three singular points, namely $A_1$, $A_2$ and
$A_3$, each of them being a $2$-fold point and the center of a quadratic
branch $P_i\in\bar\cC'$, $i=1,2,3$ \cite[Prop. 8]{ck1}. Moreover if $t$ is
a local parameter at $P_3$, a primitive representation of $P_3$ is given
by (loc. cit)
\begin{equation}\label{eq5.1}
x=t^2,\qquad y=\sum_{i=0}^{\infty} t^{2\sq+i(q-\sq+1)}\, .
\end{equation}
Since $\beta$ leaves $\cC'$ invariant, maps $A_3$ into
$A_1$, and $A_1$ to $A_2$ we obtain the following primitive 
representations for $P_1$ and $P_2$ respectively:
\begin{alignat*}{2}
x & =t^{-2\sq}+\ldots &\qquad y & =t^{-(2\sq-2)}+\ldots\\
\intertext{and}
x & =t^{2\sq-2}+\ldots &\qquad y & =t^{-2}+\ldots\, .
\end{alignat*}
Hence we have the
     \begin{lemma}\label{lemma5.1} The divisors of $x, y\in
\bar\lq(\bar\cC')$ 
are respectively:
\begin{align*}
\div(x) & = (2\sq-2)P_2+2P_3-2\sq P_1\\
\intertext{and}
\div(y) & = 2\sq P_3 - (2\sq-2)P_1-2P_2\, .
\end{align*}
In particular $x, y\in \lq(\bar\cC')$.
      \end{lemma}
Next consider the morphism 
$$
\pi=(1:x^d:y^d): \bar\cC'\to \P^2(\bar\lq)\, .
$$
Since $\pi(P)=(t^{e_P}(P):(t^{e_P}x^d)(P):(t^{e_P}y^d)(P))$, where $t$ is
a local parameter at $P$ and $e_P:=\min{(0, dv_P(x), dv_P(y))}$, $v_P$
being the valuation at $P$, from the
previous lemma we see that $\pi$ is totally ramified at $P_1$, $P_2$ and
$P_3$. Set $\cX':=\pi(\bar\cC')$.
      \begin{lemma}\label{lemma5.2}
The induced morphism $\bar\pi: \bar\cC'\to \bar\cX'$ is a $d$-sheeted
covering defined over $\lq$ which is ramified precisely at $P_1$, $P_2$
and $P_3$. Moreover, $\bar\pi$ is totally ramified at each of these
points. 
      \end{lemma}
      \begin{proof} The morphism $\bar\pi$ is defined over $\lq$ by Lemma
\ref{lemma5.1}. From (\ref{eq5.1}) we see that
$P_3'\in \bar\cX'$ over $\pi(A_3)$
has a representation of type $x'=t^{2d}$ and
$y'=(\sum_{i=0}^{\infty}t^{2\sq+i(q-\sq+1)})^d$. Putting $\tau=t^d$ we see
that $P_3'$ has a primitive representation of type
$$
x'=\tau^2\qquad y'=\tau^{2\sq}+d\tau^{2\sq+(q-\sq+1)/d}+\ldots\, ,
$$
showing that $\bar\pi$ has degree $d$. Let $P\in
\bar\cC'\setminus\{P_1,P_2,P_3\}$, hence a point over $(a:b:1)\in
\cC'$ with 
$a\neq 0$ and $b\neq 0$. Then $P$ has a primitive representation of type
$$
x=a+t\qquad y=b+b_1t+\ldots\, ,
$$ 
with $b_1\neq 0$ so that a point $P'\in \bar\cX'$ over
$\pi(a:b:1)$ has a primitive representation of type
$$
x'=a^d+dt+\ldots\qquad y'=b^d+db^{d-1}b_1t+\ldots\, .
$$
This shows that $\#\bar\pi^{-1}(\bar\pi(P))=d$ and the proof is complete.
     \end{proof}
     \begin{proof}({\em Theorem \ref{thm5.1}}) 
We show that $\bar\cX'$ above, satisfies the theorem. This curve is 
maximal over $\lq$ by Lemma \ref{lemma4.1.1},
Lemma \ref{lemma4.1.3} and \cite[Prop. 6]{lachaud}. To compute the genus
$g$ of $\bar\cX'$ we apply the Riemann-Hurwitz formula for
$\bar\pi:\bar\cC'\to \bar\cX'$ taking into 
consideration Lemma \ref{lemma5.2}. We have
$$
\sq(\sq-1)-2=d(2g-2)+3(d-1)\, ,
$$
which gives $g=\frac{1}{2}(\frac{q-\sq+1}{d}-1)$. 
     \end{proof}
     \begin{remark}\label{remark5.2} In the above notations, we have a
commutative diagram 
\begin{equation*}
\begin{CD}
\cH=\bar\cC @>\kappa>>       \cH= \bar\cC'\\
  @V{\tilde\pi}VV     @VV{\bar\pi}V\\
\bar\cX @>\kappa>>       \bar\cX'
\end{CD}
\end{equation*}
where $\cX:=\kappa^{-1}(\cX')$ and $\tilde\pi$ is the morphism induced by
$\kappa^{-1}\circ\pi\circ\kappa$. We have shown indeed that $\bar\cX'$ is
a maximal curve over $\lq\cong\fq$. A plane model for $\cX$ is giving by
$$
F(X,Y,Z)=F'(aX+Y+a^{q+1}Z, a^{q+1}X+aY+Z, X+a^{q+1}Y+aZ)=0\, ,
$$
where $F'(X_0,X_1,X_2)=0$ is a minimal equation for $\cX'$. As in the
proof of Lemma \ref{lemma4.1.3}, we see that $\cX$ is defined over $\fq$
iff $\beta(\cX')=\cX'$. For $d=3$ and using another model for $\cH$ we are
able to write down an explicit equation for $\cX$ over $\fq$ (see next
section).
     \end{remark}
To conclude this section we state some remarks on the dimension of
$\cD_{\bar\cX'}=(\sq+1)P_0$ (cf. \S\ref{s2.2}) and on Weierstrass
semigroups at some points of $\bar\cX'$. We keep the above notations.

For $i=1,2,3$, let $P_i'\in \bar\cX'$ be the unique point over $\pi(A_i)$
and let $\{Q_i\}:=\bar\pi^{-1}(P_i')$. 
Since $\pi(A_1)=A_2$, $\pi(A_2)=A_3$ and $\pi(A_3)=A_1$ and these points
do not belong to $\Pi$, then $Q_1, Q_2$ and $Q_3$ are not $\lq$-rational. 
In addition, as $p$ does not
divide $d$, the Weierstrass semigroups $H(Q_i)$ and $H(P_i')$ at $Q_i$ and
$P_i'$ respectively are related to each other as follows (see
e.g. \cite[proof of Lemma 3.4]{t})
$$
S:= H(P_i')=\{h/d: h\in H(Q_i),\ h\equiv 0\! \pmod{d}\}\, .
$$
Moreover, $H(Q_i)$ can be computed as follows \cite[Thm. 2]{g-vi},
\begin{align*}
\tilde S:=H(Q_i)= & \N\setminus\{r\sq+s+1: r+s\le \sq-2\}\\
        = & \cup_{j=1}^{\sq-2}[j\sq-(j-1),j\sq]\cup\{0,q-2\sq+2,
q-2\sq+3,\ldots\}\, .
\end{align*}
Hence from (\ref{eq2.2}) we have the
\begin{proposition}\label{prop5.1} For the curve $\bar\cX'$ above,
$$
\dim(\cD_{\bar\cX'})=1+\#\{h>0: h\in \tilde S,\ h\equiv 0\! \pmod{d},\
h\le d\sq\}\, .
$$
       \end{proposition}
       \begin{example}\label{example5.1} Let us consider the case $d=7$,
i.e. $7\mid (q-\sq+1)$ (the case $d=3$ is discussed in the next section).
By Remark \ref{remark5.1}, $\sq\equiv 3 \pmod{7}$ or $\sq\equiv 5
\pmod{7}$. The positive elements of $\tilde S$ less than or
equal to $7\sq$ are
    \begin{align*}
7\sq-6, 7\sq-5, 7\sq-4, 7\sq-3, 7\sq-2, 7\sq-1, 7\sq & \\
        6\sq-5, 6\sq-4, 6\sq-3, 6\sq-2, 6\sq-1, 6\sq & \\
                5\sq-4, 5\sq-3, 5\sq-2, 5\sq-1, 5\sq & \\
                        4\sq-3, 4\sq-2, 4\sq-1, 4\sq & \\
                                3\sq-2, 3\sq-1, 3\sq & \\
                                        2\sq-1, 2\sq & \\
                                                 \sq &\, .
    \end{align*}
\noindent So we have:

(1) If $\sq\equiv 3 \pmod{7}$, then $7\sq, 6\sq-4, 5\sq-1$ and $3\sq-2$
are the elements of $\tilde S$ which are $\le 7\sq$ and $\equiv 0
\pmod{7}$. Thus if $\sq=3$, then $\dim_{\bar\cX}=4$, and if $\sq>3$, then
$\dim(\bar\cX)=5$.

(2) If $\sq\equiv 5 \pmod{7}$, then $7\sq, 6\sq-2, 5\sq-4$ and $3\sq-1$
are the elements of $\tilde S$ which are $\le 7\sq$ and $\equiv 0
\pmod{7}$. Then $\dim(\bar\cX)=5$.

Notice that the $(\cD_{\bar\cX}, P_i')$-orders can be computed by means of 
\S\ref{s2.2}(4).
\end{example}
\section{The case $3$ divides $(q$-$\protect\sq+1)$}\label{s6}
Throughout this section we let $3$ be a divisor of $(q-\sq+1)$, or
equivalently let $\sq \equiv 2\pmod 3$. We 
keep the notations of section \ref{s4}. To construct a curve
$\bar\cX$ over $\fq$ of genus 
$(q-\sq-2)/6$, we use the equation for $\cC'=\bar\cC'$ over $\F_{\sq^3}$ 
given in Proposition \ref{prop4.2.1}. The 
construction is similar to the one given in \S\ref{s5} but we 
work in any characteristic and we shall write down an explicit equation
for $\cX$. 

Let $\cX'$ be the curve in $\P^2(\bar\fq)$ defined by
    \begin{align}\label{eq6.1}
    \begin{split}
F'(X_0,X_1,X_2)& := G(X_0,X_1,X_2)-3(X_0X_1X_2)^{\frac{\sq+1}{3}}\\
               & = X_0^{\sq}X_2+X_2^{\sq}X_1+X_1^{\sq}X_0
                   -3(X_0X_1X_2)^{\frac{\sq+1}{3}}\, .\\
    \end{split}
    \end{align}
and consider the projective morphism
$$
\pi=(X^3_0:X^3_1:X^3_2):\  \cC'=\bar\cC' \to \P^2(\bar\fq)\, .
$$ 
      \begin{lemma}\label{lemma6.1} 
    \begin{enumerate}
\item The curve $\cX'$ is geometrically irreducible.
\item $\pi(\cC')=\cX'$.
\item Let $x:=X_0/X_2, y:=X_1/X_2\in \bar\F_{\sq^3}(\cC')$. Then
\begin{align*}
\div(x) & =(\sq-1)A_2+A_3-\sq A_1\\
\intertext{and}
\div(y) & = \sq A_3-(\sq-1)A_1-A_2\, .
\end{align*}
\item The morphism $\pi$ induces a $3$-sheeted covering $\bar\pi:\cC'\to
\bar\cX'$ ramified only at $A_1$, $A_2$ and $A_3$. Moreover $\bar\pi$ is
totally ramified at each of these points.
\item The genus of $\bar\cX'$ is $\frac{1}{2}(\frac{q-\sq+1}{3}-1)$.
    \end{enumerate}
    \end{lemma}
    \begin{proof} (1) Each fundamental line $X_i=0$ ($i=0,1,2$) meets
$\cX'$ in only two points. Assume $\cX_1'$ be a component of $\cX'$.
Clearly $\cX_1'$ must meet each fundamental line in at least one point.
This implies that at least two points in the set $S:= \{A_1, A_2, A_3\}$
lie on $\cX_1'$. If we would had a further component $\cX_2'$ of $\cX'$,
then this component would also pass through two points of $S$ and one of
these points certainly be a common point of $\cX_1'$ and $\cX_2'$. Then a
point of $S$ would be a singular point of $\cX'$ which is not the case, as
a direct computation shows.

(2) From the identity
$$
a^3+b^3+c^3-3abc=(a+b+c)(\epsilon^2 a+b+\epsilon c)(\epsilon a +b +
\epsilon^2 c)\, , 
$$
where $\epsilon$ is a primitive $3$-th root of unity, we
have that 
$$
F'(X_0^3, X_1^3, X_2^3)=G(X_0,X_1,X_2)G(\epsilon X_0,\epsilon X_1,X_2)
G(X_0, X_1, \epsilon X_2)\, .
$$
This and (1) implies (2).

(3) Similar to the proof of Lemma \ref{lemma5.1}.

(4) Similar to the proof of Lemma \ref{lemma5.2}.

(5) Follows from (4) and the Riemann-Hurwitz formula applied to $\bar\pi$.
     \end{proof}
Consider the following commutative diagram 
\begin{equation*}
\begin{CD}
\cH @>\kappa>>        \cC'=\bar\cC'\\
  @V{\tilde\pi}VV     @VV{\bar\pi}V\\
\bar\cX @>\kappa>>       \bar\cX'
\end{CD}
\end{equation*}
where $\cX:= \kappa^{-1}(\cX')$ and $\tilde \pi$ is the morphism induced
by $\kappa^{-1}\circ\pi\circ\kappa$. Now we can state the main result of
this section.
          \begin{theorem}\label{thm6.1} The curve $\bar\cX$ above is a 
maximal curve over $\fq$ of genus $(q-\sq-2)/6$. A plane model over $\fq$
for $\cX$ is given by 
\begin{equation}\label{eq6.2}
F(X,Y,Z)=cF'(aX+Y+a^{\sq+1}Z, a^{q+1}X+aY+Z,X+a^{q+1}Y+aZ)=0\, ,
\end{equation}
where $a\in \F_{\sq^3}$ satisfies Lemma \ref{lemma4.2.1}, $c\in \F_{q^3}$
such that $c^{\sq-1}=a$, and $F'(X_0, X_1,X_2)$ is the polynomial in
(\ref{eq6.1}).
          \end{theorem}
\begin{proof} The statement on the genus follows from Lemma 
\ref{lemma6.1}(5) and the diagram above. To see that $\bar\cX$ is
maximal, by \cite[Proposition 6]{lachaud}, it is enough to show that
$\cX$ is defined over $\fq$. Let $G'$ be a
minimal equation of $\cX$ so that
$$
G'(X,Y,Z):=F'(aX+Y+a^{q+1}Z,a^{q+1}X+aY+Z,X+a^{q+1}Y+aZ)=0\, .
$$
Note that $F'(X_0,X_1,X_2)\in \fq [X_0,X_1,X_2]$  
and that $F'(X_0,X_1,X_2)=F'(X_2,X_0,X_1)$, i.e. $\beta(\cX')=\cX')$. 
Using these facts we obtain
\begin{align*}
G'(X,Y,Z)^q  
= & 
F'(a^qX^q+Y^q+a^{q^2+q}Z^q,a^{q^2+q}X^q+a^qY^q+Z^q,X^q+a^{q^2+q}Y^q+a^qZ^q)\\
= &
F'(X^q+a^{q^2+q}Y^q+a^qZ^q,a^{q}X^q+Y^q+a^{q^2+q}Z^q,
a^{q^2+q}X^q+a^{q}Y^q+Z^q)\, .
\end{align*}
Now from $a^{q+\sq+1}=1$ (cf. proof of Lemma \ref{lemma4.2.1}) we get
$a^{q^2+q+1}=1$ and since 
$F'(X_0,X_1,X_2)$ is a homogeneous polynomial of degree $\sq+1$ we have
that 
$$
G'(X,Y,Z)^q=a^{-(\sq+1)}G'(X^q,Y^q,Z^q)\, .
$$
Then $\cX$ can also be defined by 
$$
F(X,Y,Z):= cG'(X,Y,Z)\qquad \text{with\quad $c^{\sq-1}=a$}\, ,
$$
(so that $c\in F_{q^3}$ as $a^{q^2+q+1}=1$), and since
$F(X,Y,Z)^q=F(X^q,Y^q,Z^q)$ we are done.
\end{proof}
Next we investigate the linear system $\cD_{\bar\cX}:=|(\sq+1)P_0|$
(cf. \S\ref{s2.2}). For $i=1,2,3$, let $P_i'\in \bar\cX'$ be the unique
point over $\pi(A_i)$ and let $P_i:=\kappa^{-1}(P_i')\in
\bar\cX$. We notice that $P_i\not\in
\bar\cX(\fq)$ because $\pi(A_i)\in \{A_1, A_2,A_3\}$ and the coordinates
of $\kappa^{-1}(A_i)$ do not belong to $\fq$.
\begin{lemma}\label{lemma6.2} The first two positive Weierstrass non-gaps
at $P_i$ above ($i=1,2,3$) are $(2\sq-1)/3$ and $\sq$.
\end{lemma}
\begin{proof} Is a particular case of the discussion after Remark
\ref{remark5.2}.
\end{proof}
    \begin{proposition}\label{prop6.1} 
Let $\bar\cX$ be the maximal curve over $\fq$ of Theorem \ref{thm6.1}. 
Then
\begin{enumerate}
\item $\dim(\cD_{\bar\cX})=3$.
\item The $(\cD_{\bar\cX},P)$-orders are:

(i) $0,1,2$ and $\sq+1$ if $P\in \bar\cX(\fq$;

(ii) $0,1, (\sq+1)/3$ and $\sq$ if $P\in \{P_1, P_2,P_3\}$;

(iii) $0,1,2$ and $\sq$ if $P\in
\bar\cX\setminus(\bar\cX(\fq)\cup\{P_1,P_2,P_3\})$.
\end{enumerate}
    \end{proposition}
     \begin{proof} Set $\cD:=\cD_{\bar\cX}$.

(1) It follows from Lemma \ref{lemma3.1}(1) (or from
Prop. \ref{prop5.1} with $d=3$).

(2) The $\cD$-orders are 0, 1, 2 and $\sq$ by Lemma \ref{lemma3.1}(2). 
From Lemma \ref{lemma6.2} and \S\ref{s2.2}(1)(4) for each $P_i$ 
($i=1,2,3$) the $(\cD,P_i)$-orders are $0,1, (\sq+1)/3$ and $\sq$. Hence
from \S\ref{s2.1}(3), it follows that $v_{P_i}(R^{\cD})=(\sq-5)/3$. Now, 
using $g=(q-\sq-2)/6$ and the maximality of $\bar\cX$ we obtain
$$
\deg(R^{\cD})-3(\sq-5)/3=\#\bar\cX(\fq)\, .
$$
Then (2) follows by taking into consideration that $j_3(P)=\sq+1$ at each
$\fq$-rational point (cf. \S\ref{s2.2}(3)).
\end{proof}
\begin{remark}\label{remark6.1} We can describe $\cD_{\bar\cX}$ via
hyperplane sections on a certain curve. Indeed, consider the morphism
$$
\varphi=(1:x^3:y^3:xy):\ \cC'\to \P^3(\bar\F_{\sq^3})
$$
with $x=X_0/X_2, y=X_1/X_2\in \bar\F_{\sq^3}(\cC')$ satisfying
$x^{\sq}+y+xy^{\sq}=0$. Let $\cY:=\varphi(\cC')$ and let ${\rm Pr}$
denotes the projection of  
$\P^3(\bar\fq)\setminus\{(0:0:0:1)$ from the point $(0:0:0:1)$ on the
plane $\P^2(\fq)\cong \{(x:y:z:0): (x:y:z)\in \P^2(\bar\fq)\}$. Then
arguing as in Lemma \ref{lemma5.2} one can show that
    \begin{enumerate}
\item[(i)] \quad $\cX'={\rm Pr}(\cY)$; 
\item[(ii)]\quad  $\cD=(\kappa^{-1}\circ {\rm Pr})^{*}(\cD_{\bar\cX})$,
    \end{enumerate}
where $\cD$ cuts $\cY$ by hyperplanes of $\P^3(\bar\F_{\sq^3})$.
\end{remark}
So far we have the following maximal curves $\cX$ over $\fq$ such that
$\dim(\cD_{\cX})=3$:
\begin{enumerate}
\item The non-singular model of $y^{\sq}+y=x^{(\sq+1)/2}$, $q$ odd,
whose genus is $(\sq-1)^2/4$ (cf. Remark \ref{remark3.2}(1));
\item The non-singular model of $\sum_{i=0}^{t}y^{\sq/2^i}=x^{\sq+1}$,
$q=2^t$, whose genus is \\
$\sq(\sq-2)/4$ (cf. Remark \ref{remark3.2}(2);
\item A maximal curve of genus $\sq(\sq-1)/6$, $\sq\equiv 0\pmod{3}$, (cf.
Remark \ref{remark3.4}).
\item A maximal curve of genus $(q-\sq-2)/6$ (cf. Theorem \ref{thm6.1}).
\end{enumerate}
A natural question then is the following:
\begin{question*}
Are all the maximal curves $\cX$ over $\fq$ (up to $\fq$-isomorphism)
with $\dim(\cD_{\cX})=3$ listed above ?
\end{question*}
\begin{remark}\label{remark6.2} Connection between maximal curves and
non-classical curves have been noticed in \cite[Prop. 1.7]{fgt}. More
precisely, a maximal curve over $\fq$ of genus $g$ is
non-classical provided that $g\ge \sq-1$. Hence the examples obtained here
are non-classical for $d\le (q-\sq+1)/(2\sq-1)$. Most of the known
examples of non-classical curves are Artin-Schreier extensions of rational
function fields, with genus a multiple of the characteristic \cite{sch},
\cite{g-vi}. Thus the maximal curves obtained here are in fact new
examples of non-classical curves.
\end{remark}

\end{document}